\documentclass[10pt]{amsart}

\usepackage[alphabetic]{amsrefs}
\usepackage[all]{xy}
\usepackage{amssymb,amsmath, latexsym}

\newtheorem{lemma}{Lemma}[section]
\theoremstyle{definition}
\newtheorem{question}[lemma]{Question}
\newtheorem{definition}[lemma]{\sl Definition}
\newtheorem{notation}[lemma]{\sl Notation}
\newtheorem{definition-notation}[lemma]{\sl Definition/Notation}

\newtheorem{proposition}[lemma]{Proposition}
\newtheorem{theorem}[lemma]{Theorem}

\newtheorem{corollary}[lemma]{Corollary}

\theoremstyle{remark}
\newtheorem{remark}[lemma]{Remark}

\newenvironment{pf}{\noindent{\textbf{Proof.}}}{\hfill $\square$\medskip}

\def\ra{\rightarrow}
\def\w{\wedge}
\def\bz{\mathbb{Z}}
\def\bn{\mathbb{N}}
\def\bf{\mathbb{F}}
\def\<{\langle}
\def\>{\rangle}
\def\loc{{\sf loc}}
\def\lim{{\sf lim}}
\def\th{{\sf th}}
\def\S{{\sf S}}
\def\T{{\sf T}}
\def\C{{\sf C}}
\def\FF{\mathcal{F}}
\def\CC{\mathcal{C}}
\def\BA{\mathsf{BA}}
\def\BL{\mathsf{BL}}
\def\DL{\mathsf{DL}}
\def\0{\<0\>}
\def\SS{\mathcal{S}}
\def\L{\mathcal{L}}
\def\A{\mathcal{A}}
\def\TC{{\sf TC}}
\def\TCi{{\sf TC1_n}}
\def\TCij{{\sf TC1_i}}
\def\TCii{{\sf TC2_n}}
\def\TCiijj{{\sf TC2_{n-1}}}
\def\TCiii{{\sf TC3_n}}

\def\GSC{{\sf GSC}}
\def\SDGSC{{\sf SDGSC}}
\def\iff{\Leftrightarrow}
\def\LTC{{\sf LTC}}
\def\LTCi{{\sf LTC1_n}}
\def\LTCij{{\sf LTC1_i}}
\def\LTCii{{\sf LTC2_n}}
\def\LTCiij{{\sf LTC2_i}}
\def\LTCiijj{{\sf LTC2_{n-1}}}
\def\LTCiii{{\sf LTC3_n}}
\def\LTCiiij{{\sf LTC3_i}}
\def\Ln{L_n}
\def\Lnf{L_n^f}
\def\lnf{{\it l}_n^f}
\def\ln{{\it l}_n}

\def\l{{\it l}}
\def\c{{\it c}}
\def\HH{\mathcal{H}}
\def\LLn{\mathcal{L}_n}
\def\KKn{\mathcal{K}_n}
\def\LL{\mathcal{L}}

\newcommand{\boldone}{{\rm 1\hspace*{-0.4ex}%
\rule{0.1ex}{1.52ex}\hspace*{0.2ex}}}

\def\one{\boldone}
\def\hfp{{H\bf_p}}
\def\ua{{\uparrow}}

\def\bp{{BP}}
\def\supp{{\sf supp}}

\begin{document}

\title[Variations of the telescope conjecture]{Variations of the telescope conjecture and Bousfield lattices for localized categories of spectra}
\author{F. Luke Wolcott}             
\email{luke.wolcott@gmail.com}
\address{Department of Mathematics,
         Lawrence University,
         Appleton, WI, 54911,
         USA}
%
%
%\address{Department of Mathematics,
%         University of Western Ontario,
%         London, Ontario N6A~5B7,
%          Canada}

\date{\today}

\begin{abstract}
We investigate several versions of the telescope conjecture on localized categories of spectra, and implications between them.  Generalizing the ``finite localization" construction, we show that on such categories, localizing away from a set of strongly dualizable objects is smashing.  We classify all smashing localizations on the harmonic category, $H\bf_p$-local category and $I$-local category, where $I$ is the Brown-Comenetz dual of the sphere spectrum; all are localizations away from strongly dualizable objects, although these categories have no nonzero compact objects.  The Bousfield lattices of the harmonic, $E(n)$-local, $K(n)$-local, $H\bf_p$-local and $I$-local categories are described, along with some lattice maps between them.  One consequence is that in none of these categories is there a nonzero object that squares to zero.  Another is that the $H\bf_p$-local category has localizing subcategories that are not Bousfield classes.
\end{abstract}

\maketitle

\section{Introduction}
\vspace{22pt}

The telescope conjecture, first stated by Ravenel~\cite[Conj.~10.5]{[Rav84]}, is a claim about two classes of localization functors in the $p$-local stable homotopy category of spectra.  First, one can localize away from a finite type $n+1$ spectrum $F(n+1)$; the acyclics are the smallest localizing subcategory containing $F(n+1)$, and we denote this functor $L_{n}^f$.  Second, one can localize at the wedge of the first $n+1$ Morava $K$-theories $K(0)\vee\cdots \vee K(n)$; the acyclics are all spectra that smash with $K(0)\vee\cdots\vee K(n)$ to zero and this is denoted $L_{n}$.  Both these localizations are smashing, i.e.~they commute with coproducts.  The telescope conjecture ($\TC_n$), basically, claims that $\Lnf$ and $\Ln$ are isomorphic.  In fact, here we consider three slightly different versions, $\TCi$, $\TCii$, and $\TCiii$, of the telescope conjecture.  In Section~\ref{sect-TCs} we articulate them carefully and show implications between them.

The conjecture is known to hold for $n=0$~\cite[p.~79]{[Rav92]}, and for $n=1$ when $p=2$~\cite{[Mah82]} and $p>2$~\cite{[Mil81]}.  A valiant but unsuccessful effort at a counterexample, for $n\geq 2$, was undertaken by Mahowald, Ravenel, and Shick, as outlined in~\cite{[MRS]}.  Since then little progress has been made, and the original conjecture remains open.

A generalization of the telescope conjecture can be stated for spectra, as well as other triangulated categories.  Localization away from a finite spectrum, i.e.~a compact object of the category, always yields a smashing localization functor (see e.g.~\cite[Prop.~2.9]{[Bou79a]} or~\cite{[Miller]} or~\cite[Thm.~3.3.3]{[HPS]}). The Generalized Smashing Conjecture ($\GSC$) is that every smashing localization arises in this way.  If true, then every smashing localization is determined by its compact acyclics; if the $\GSC$ holds in spectra, then so must the $\TC_n$ for all $n$.  

The $\GSC$, essentially stated for spectra decades ago by Bousfield \cite[Conj.~3.4]{[Bou79b]}, has been formulated in many other triangulated categories, in many cases labeled as the telescope conjecture, and in many cases proven to hold.  Neeman~\cite{[Neeman92]} made the conjecture for the derived category $D(R)$ of a commutative ring $R$, and showed it holds when the ring is Noetherian.  See also~\cite[Thm.~6.3.7]{[HPS]} or~\cite{[KrSt10]} for a generalization.  On the other hand, Keller~\cite{[keller]} gave an example of a non-Noetherian ring for which the $\GSC$ fails.  Benson, Iyengar, and Krause have shown that the $\GSC$ holds in a stratified category~\cite{[BIKstcat]}, such as the stable module category of a finite group~\cite{[BIKstcat2]}.  Balmer and Favi~\cite{[BF]} show that in a tensor triangulated category with a good notion of support, the $\GSC$ is a ``local" question.

It is worth noting that there are further variations of the \GSC, which we won't consider here.  Krause~\cite{[krause-sstc]} formulated a variation of the $\GSC$, in terms of subcategories generated by sets of maps, that makes sense (and holds) for any compactly generated triangulated category.  Krause and Solberg give a variation for stable module categories, stated in terms of cotorsion pairs~\cite{[KS03]}.  See also~\cites{[krause-cqsl],[HST08],[Bruning07],[Stov10]}.

To date, Keller's ring yields the only category where the $\GSC$ is known to fail.  In this paper we give several more examples.  Incidentally, each is a well generated triangulated category that is not compactly generated.

One of our main results is the following.  We weaken the assumptions for ``finite localization", and show that in many categories, localization away from any set of strongly dualizable objects yields a smashing localization.  (Recall that an object $X$ is strongly dualizable if $F(X,Y) \cong F(X,\one)\w Y$ for all $Y$, where $\one$ is the tensor unit and $F(-,-)$ the function object bifunctor.)  Let $\loc(X)$ denote the smallest localizing subcategory containing $X$.  We prove the following as Theorem~\ref{sdsl}.\\

\noindent \textbf{Theorem A:} Let $\T$ be a well generated tensor triangulated category such that  $\loc(\one) = \T$.  Let $A=\{B_\alpha\}$ be a (possibly infinite) set of strongly dualizable objects.  Then there exists a smashing localization functor $L:\T\ra \T$ with Ker~$L = \loc(A)$. \\

Thus we are led to conjecture the following.\\

\noindent \textbf{Strongly Dualizable Generalized Smashing Conjecture ($\SDGSC$):} Every smashing localization is localization away from a set of strongly dualizable objects.\\

We give several examples of categories where the $\GSC$ fails, but the $\SDGSC$ holds.  In fact, we consider a topological setting, where one can also formulate a version (or versions, rather) of the original telescope conjecture.  

Specifically, we consider localized categories of spectra.  Let $\SS$ be the $p$-local stable homotopy category, and let $\w$ denote the smash, i.e. tensor, product.  Take $Z$ an object of $\SS$, and let $L=L_Z:\SS\ra \SS$ be the localization functor that annihilates $Z_*$-acyclic objects.  The full subcategory of $L$-local objects, that is, objects $X$ for which $X {\ra} LX$ is an equivalence, has a tensor triangulated structure induced by that of $\SS$.  Let $\LL$ denote this category; the triangles are the same as in $\SS$, the coproduct is $X\coprod_\LL Y = L(X\coprod Y)$ and the tensor is $X\w_\LL Y = L(X\w Y)$.  

In Definition~\ref{defn-l}, we define localization functors $\lnf$ and $\ln$ on $\LL$ that are localized versions of $\Lnf$ and $\Ln$.  The localized telescope conjecture ($\LTC$), basically, is that $\lnf$ and $\ln$ are isomorphic.  In fact, we give three versions of the $\LTC$, and in Theorems~\ref{LTCs} and~\ref{LTCstwo} establish implications between them.  Then, examining specific examples of localized categories of spectra, we conclude the following in Theorems~\ref{LTCs-harm}, \ref{LTCs-kn}, \ref{LTCs-hfp}, \ref{i-loc-tc} and \ref{LTCs-bp} and Corollary~\ref{LTCs-en}.  \\

\noindent \textbf{Theorem B:} All versions of the localized telescope conjecture, $\LTCij$, $\LTCiij$, and $\LTCiiij$ hold for all $i\geq 0$, in the $\bigvee_{n\geq 0}K(n)$-local (i.e. harmonic), $K(n)$-local, $\hfp$-local, $BP$-local, and $I$-local categories, where $I$ is the Brown-Comenetz dual of the sphere spectrum.\\

In order to consider the $\GSC$ and $\SDGSC$ in $\LL$, we must classify the smashing localizations on $\LL$.  We are able to do this in several examples.\\

\noindent \textbf{Theorem C:} In the harmonic category, the $\GSC$ fails but the $\SDGSC$ holds.  Likewise in the $\hfp$-local and $I$-local categories.  In the $BP$-local category the $\GSC$ fails but the $\SDGSC$ is open.  In the $E(n)$-local and $K(n)$-local categories the $\GSC$ and $\SDGSC$ both hold.\\

\begin{pf} This is Theorems~\ref{harm-sdgsc} and \ref{LTCs-kn}, Propositions~\ref{gsc-hfp} and \ref{gsc-bp}, and Corollaries~\ref{LTCs-en} and \ref{gsc-i}. \end{pf}

One novelty in our approach is our use of Bousfield lattice arguments.  Given an object $X$ in a tensor triangulated category $\T$, the Bousfield class of $X$ is $\<X\> = \{W\;|\; W\w X =0\}$.  It is now known~\cite{[IK]} that every well generated tensor triangulated category has a set of Bousfield classes.  This set has the structure of a lattice, and is called the Bousfield lattice of $\T$.  One can now attempt to calculate the Bousfield lattices of categories of localized spectra.  Furthermore, every smashing localization yields a pair of so-called complemented Bousfield classes.  Information about the Bousfield lattice of a category gives information about its complemented classes, which gives information about the smashing localization functors on the category.

Moreover, the first version of the telescope conjecture $\TCi$ is that two spectra $T(n)$ and $K(n)$ have the same Bousfield class.  In the localized version this becomes ($\LTCi$) the claim that $\<LT(n)\>=\<LK(n)\>$ in the Bousfield lattice of $\LL$.  One is thus led to investigating Bousfield lattices of localized spectra.

Corollary~\ref{bl-size} gives an upper bound, $2^{2^{\aleph_0}}$, on the cardinality of such lattices.  Jon Beardsley has calculated the Bousfield lattice of the harmonic category to be isomorphic to the power set of $\bn$; we give this calculation in Proposition~\ref{jb}.  In Corollary~\ref{ln-bl} and Proposition~\ref{realize} we show that one can realize this lattice as an inverse limit of the Bousfield lattices of $E(n)$-local categories, as $n$ ranges over $\bn$.  Then in Corollary~\ref{bl-kn} and Propositions~\ref{bl-hfp} and ~\ref{bl-i}, we show that the $K(n)$-local, $\hfp$-local, and $I$-local categories all have two-element Bousfield lattices.  In Proposition~\ref{bl-bp} we give a lower bound, $2^{\aleph_0}$, on the cardinality of the Bousfield lattice of the $BP$-local category.

%As an immediate object-level application of these Bousfield lattice calculations, we have the following.  Call an object $X$ {\it square-zero} if it is nonzero, but $X\w X=0$.  For example, $I\w I=0$ in $\SS$.
%
%\begin{proposition} There are no square-zero objects in the harmonic, $E(n)$-, $K(n)$-, $\hfp$-, or $I$-local categories.
%\end{proposition}
%
%\begin{pf} In Corollary 2.8 of~\cite{[preprint]}, we show that in a well generated tensor triangulated category, there are no square-zero objects if and only if $\BA=\DL=\BL$.  The claim then follows from Corollaries~\ref{ln-bl} and~\ref{bl-kn} and Propositions~\ref{bl-hfp} and~\ref{bl-i}.
%\end{pf}

One immediate object-level application of these Bousfield lattice calculations is the following.  Call an object $X$ square-zero if it is nonzero, but $X\w X=0$.  Then Proposition~\ref{no-sz} shows that there are no square-zero objects in the harmonic, $E(n)$-, $K(n)$-, $\hfp$-, or $I$-local categories.

We are also able to answer the analogue of a conjecture by Hovey and Palmieri, originally stated for the stable homotopy category.  Conjecture 9.1 in~\cite{[HP]} is that every localizing subcategory is a Bousfield lattice.  Proposition~\ref{loc-not-bc} demonstrates that this fails in the $\hfp$-local category, by giving two localizing subcategories that are not Bousfield classes.

Section~\ref{sect-prelims} establishes the categorical setting, and provides background on localization, Bousfield lattices, and stable homotopy theory.  Section~\ref{sect-TCs} defines the various versions of the telescope conjecture, for spectra and for localized spectra, and establishes implications among them.  The remainder of the paper is devoted to examining specific examples: the harmonic category (Section~\ref{sect-harm}), the $E(n)$-local and $K(n)$-local categories (Section~\ref{sect-enkn}), and the $\hfp$-local, $I$-local, $BP$-local, and $F(n)$-local categories (Section~\ref{sect-other}).  All results are new unless cited.  Most of the results on the $E(n)$-local and $K(n)$-local categories in Section~\ref{sect-enkn} follow in a straightforward way from Hovey and Strickland's work in~\cite{[HovS]}, and are included for completeness.  \\

We would like to thank Dan Christensen for extensive discussions and suggestions, and Jon Beardsley for Proposition~\ref{jb}.

\vspace{22pt}
\section{Preliminaries} \label{sect-prelims}
\vspace{22pt}

\subsection{Categorical setting} We start with the notion of a tensor triangulated category $\C$; i.e.~ a triangulated category with set-indexed coproducts and a closed symmetric monoidal structure compatible with the triangulation~\cite[App.A]{[HPS]}.  Let $\Sigma: \C\ra\C$ denote the shift, $[X,Y]$ the morphisms from $X$ to $Y$, and $[X, Y]_n = [\Sigma^n X, Y]$ for any $n\in\bz$.

Let $-\w-$ denote the smash (tensor) product, $\one$ denote the unit, and $F(-,-)$ denote the function object bifunctor; $F(X,-)$ is the right adjoint to $X\w-$.  Recall that an object $X$ in $\C$ is said to be {\it strongly dualizable} if the natural map $DX\w Y\ra F(X,Y)$ is an isomorphism for all $Y$, where $DX=F(X,\one)$ is the Spanier-Whitehead dual.  Since $F(\one, X)\cong X$ for all $X$, the map $F(\one,\one)\w Y\ra F(\one, Y)$ is an equivalence and $\one$ is always strongly dualizable.

For a regular cardinal $\alpha$, we say an object $X$ is {\it $\alpha$-small} if every morphism $X\ra \coprod_{i\in I} Y_i$ factors through $\coprod_{i\in J}Y_i$ for some $J\subseteq I$ with $|J|<\alpha$.  If $X$ is $\aleph_0$-small we say $X$ is {\it compact} (\cite{[HPS]} calls this {\it small}); this is equivalent to the condition that the natural map $\oplus_{i\in K}[X,Z_i]\ra [X,\coprod_{i\in K} Z_i]$ is an isomorphism for any set-indexed coproduct $\coprod_{i\in K} Z_i$.  We say $\C$ is {\it $\alpha$-well generated} if it has a set of perfect generators~\cite[Sect.~5.1]{[krause-loc-survey]} which are $\alpha$-small, and $\C$ is {\it well generated} if it is $\alpha$-well generated for some $\alpha$.  See~\cite{[krause-loc-survey]} for more details.

A {\it localizing} subcategory is a triangulated subcategory of $\C$ that is closed under retracts and coproducts; a {\it thick} subcategory is a triangulated subcategory that is closed under retracts.  Given an object or set of objects $X$, let $\loc(X)$ (resp. $\th(X)$) denote the smallest localizing (resp. thick) subcategory containing $X$.  We say that $\loc(X)$ is {\it generated by }$X$.

%\begin{definition-notation}  Throughout this paper let $\T$ be a well generated tensor triangulated category such that $\one$ is strongly dualizable and $\loc(\one) = \T$. \label{cat}
%\end{definition-notation}

\begin{definition-notation}  Throughout this paper let $\T$ be a well generated tensor triangulated category such that  $\loc(\one) = \T$. \label{cat}
\end{definition-notation}

In the language of~\cite{[HPS]}, such a $\T$ is almost a ``monogenic stable homotopy category", except that we do not insist that the unit $\one$ is compact.

In practice, in this paper $\T$ will always be either the $p$-local stable homotopy category of spectra $\SS$ or the category $\LL_Z$ of $L_Z$-local objects derived from a localization functor $L_Z:\SS\ra \SS$.  The former satisfies Definition
~\ref{cat} by~\cite[Ex.~1.2.3(a)]{[HPS]}, and the latter by Theorem~\ref{local-cat} and Lemma~\ref{local-cat-blah} below.\\

\subsection{Background on localization}  Recall that a {\it localization functor} (or simply {\it localization}) on a tensor triangulated category $\C$ is an exact functor $L:\C\ra \C$, along with a natural transformation $\eta: 1\ra L$ such that $L\eta$ is an equivalence and $L\eta=\eta L$.  We call Ker~$L$ the {\it $L$-acyclics}.  It follows that there is an exact functor $C:\C\ra \C$, called {\it colocalization}, such that every $X$ in $\C$ fits into an exact triangle $CX\ra X \ra LX$, with $CX$ $L$-acyclic.  An object $Y$ is {\it $L$-local} if it is in the essential image of $L$, and this is equivalent to satisfying $[Z,Y] = 0$ for all $L$-acyclic $Z$.  See~\cite[Ch.~3]{[HPS]} or~\cite{[krause-loc-survey]} for further background.

We also recall two special types of localizations.  A localization $L:\C\ra \C$ is said to be {\it smashing} if $L$ preserves coproducts, equivalently if $LX\cong L\one\w X$ for all $X$.  

Given a set $A$ of objects of $\C$, we say that a localization functor $L:\C\ra\C$ is {\it localization away from $A$} if the $L$-acyclics are precisely $\loc(A)$.  If such a localization exists, we also say it is {\it generated by $A$}.  When $\C=\SS$, it is well known (e.g.~\cites{[Miller], [MS]}) that localization away from a set of compact objects exists, and yields a smashing localization functor.  As mentioned in the introduction, this result has been generalized to other categories as well (e.g.~\cite[Thm.~3.3.3]{[HPS]}, \cite[Thm.~4.1]{[BF]}).  We present a further generalization in Theorem~\ref{sdsl}.

In this paper we will restrict our attention to homological localizations, which we now describe.  Given an object $Z$ in a tensor triangulated category $\C$, the {\it Bousfield class} of $Z$ is defined to be
$$\<Z\> = \{W\in \C\;|\; W\w Z=0\}.$$

Extending a classical result of Bousfield's for $\SS$, Iyengar and Krause recently showed~\cite[Prop.~2.1]{[IK]} that for every object $Z$ in a well generated tensor triangulated category $\C$, there is a localization functor $L_Z:\C\ra \C$ with $L_Z$-acyclics precisely $\<Z\>$.  We call such an $L_Z$ {\it homological localization at $Z$}.

\begin{notation} \label{loc-notation} Let $\T$ be as in Definition~\ref{cat}, with tensor unit $\one$.  For an object $Z$ in $\T$, let $L_Z:\T\ra \T$ be homological localization at $Z$.  Let $\LL_Z$ denote the category of $L_Z$-local objects, the essential image of $L_Z$.\\
\end{notation}

\begin{theorem}~\cite[3.5.1,3.5.2]{[HPS]} \label{local-cat} Let $L=L_Z:\T\ra \T$ be a localization, and $\LL_Z$ the category of $L_Z$-local objects.  Then $\LL_Z$ has a natural structure as a tensor triangulated category, generated by $L_Z\one$, which is the unit.  Considered as a functor from $\T$ to $\LL_Z$, $L$ preserves triangles, the tensor product and its unit, coproducts, and strong dualizability.  Furthermore, $L$ preserves compactness if and only if $L$ is a smashing localization.
\end{theorem}

Explicitly, for $L$-local objects $X$, $X_i$ and $Y$, in $\LL$ we have $\coprod_{\LL} X_i = L(\coprod_{\T} X_i)$ and $X \w_\LL Y = L(X\w_\T Y)$ and $F_\LL(X,Y) = F(X,Y)$.  Note that $L_Z\one$ is strongly dualizable but may not be compact in $\LL_Z$.

\begin{lemma} \label{local-cat-blah} The category $\LL_Z$ is well generated.
\end{lemma}

\begin{pf} By Proposition 2.1 of~\cite{[IK]}, the $L_Z$-acyclics $\<Z\>$ form a well generated localizing subcategory of $\T$.  Then by~\cite[Thm.~7.2.1]{[krause-loc-survey]}, the Verdier quotient $\T/\<Z\>$, which is equivalent to the local category $\LL_Z$, is well generated.  \end{pf}

We conclude this subsection with a lemma with four useful well-known facts.  Recall that a {\it ring object} in a tensor triangulated category is an object $R$ with an associative multiplication map $\mu: R\w R\ra R$ and a unit $\iota: \one\ra R$, making the evident diagrams commute.  If $R$ is a ring object, then an {\it $R$-module object} is an object $M$ with a map $m: R\w M\ra M$ along with evident commutative diagrams.  Note that $R\w X$ is an $R$-module object for every $X$.  A {\it skew field object} is a ring object such that every $R$-module object is {\it free}, i.e.~isomorphic to a coproduct of suspensions of $R$~\cite[Def.~3.7.1]{[HPS]}.

\begin{lemma} \label{useful} Let $\C$ be a tensor triangulated category with $\loc(\one)=\C$, and $L:\C\ra \C$ a localization.
\begin{enumerate} 
     \item Every localizing subcategory $\S$ of $\C$ is tensor-closed; that is, given $X\in \S$ and $Y\in \C$, $X\w Y\in \S$.
     \item	For all $X$ and $Y$ in $\C$, $L(X\w Y) = L(LX\w LY)$.
     \item Considered as a functor from $\C$ to $\LL$, $L$ also preserves ring objects and  module objects. 
		 \item If $R$ is a ring object and $M$ is an $R$-module object (in particular, if $M=R$), then $M$ is $R$-local. 
\end{enumerate}
\end{lemma} 

\begin{pf} For $(1)$, note that $Y\in \loc(\one) = \C$, so $X\w Y\in \loc( X\w \one)=\loc(X)\subseteq \S$.

For $(2)$, consider the exact triangle $X\w CY\ra X\w Y\ra X\w LY$.  Since $CY$ is $L$-acyclic and these form a localizing subcategory, $L(X\w CY)=0$, so $L(X\w Y) = L(X\w LY)$.  Using the same reasoning with the triangle $CX\w LY\ra X\w LY\ra LX\w LY$, the result follows.

 If $R \in \C$, is a ring object, then $L(\mu): L(R\w R) = L(LR\w LR) = LR\w_{\LL} LR\ra LR$, and all the localized diagrams commute.  Similarly for module objects.  

Part (4) is~\cite[Prop.~1.17(a)]{[Rav84]}. \end{pf}

\subsection{Background on Bousfield lattices} Every well generated tensor triangulated category, and hence every localized category of spectra, has a set (rather than a proper class) of Bousfield classes~\cite[Thm.~3.1]{[IK]}.  This was also recently shown for every tensor triangulated category with a combinatorial model~\cite{[CGR]}.  This set is called the {\it Bousfield lattice} $\BL(\T)$ and has a lattice structure which we now recall.  Refer to~\cites{[HP],[preprint]} for more details.

The partial ordering is given by reverse inclusion: we say $\<X\>\leq \<Y\>$ when $W\w Y=0 \implies W\w X=0$.  It is also helpful to remember that, unwinding definitions, $\<X\>\leq \<Y\>$ precisely when every $L_X$-local object is also $L_Y$-local.  Clearly $\0$ is the minimum and $\<\one\>$ is the maximum class.  The join of any set of classes is $\bigvee_{i\in I} \<X_i\> = \<\coprod_{i\in I} X_i\>$, and the meet is defined to be the join of all lower bounds.  

The smash product induces an operation on Bousfield classes, where $\<X\>\w \<Y\> = \<X\w Y\>$.  This is a lower bound, but in general not the meet.  However, if we restrict to the subposet $\DL = \{\<W\>\;|\; \<W\w W\> = \<W\>\}$, then the meet and smash agree.  Since coproducts distribute across the smash product, $\DL$ is a distributive lattice.

We say a class $\<X\>$ is {\it complemented} if there exists a class $\<X^c\>$ such that $\<X\>\w \<X^c\>=\0$ and $\<X\>\vee\<X^c\> = \<\one\>$.  The collection of complemented classes is denoted $\BA$.  For example, every smashing localization $L:\T\ra \T$ gives a pair of complemented classes, namely $\<C\one\>$ and $\<L\one\>$.  Because every complemented class is also in $\DL$, $\BA$ is a Boolean algebra.

\begin{proposition} \label{lattice-map} Let $\T$ be as in Definition~\ref{cat}, and $L_Z:\T\ra \T$ a localization functor as in Notation~\ref{loc-notation}.  Then $L_Z$ induces a well-defined order-preserving map of lattices $\BL(\T)\ra \BL(\LL_Z)$, where $\<X\>\mapsto \<L_ZX\>$.  This map is surjective, and sends $\DL(\T)$ onto $\DL(\LL_Z)$ and $\BA(\T)$ onto $\BA(\LL_Z)$.
\end{proposition}

\begin{pf} Most of this was proved in Lemma $3.1$ of~\cite{[preprint]}.  For $\<X\>\in \DL(\T)$, using Lemma~\ref{ugly} we get
$$\<LX\> = \<L(X\w X)\> = \<L(LX\w LX)\> = \<LX \w_\LL LX\>.$$  Likewise, one can check that for $\<X\>\in \BA(\T)$, the class $\<LX\>\in \BL(\LL)$ is complemented by $\<L(X^c)\>$, keeping in mind that $\<L\one\>$ is the top class in $\BL(\LL)$. \end{pf}

\begin{corollary} \label{bl-size} For any $Z\in \SS$, we have $|\BL(\LL_Z)|\leq 2^{2^{\aleph_0}}$.
\end{corollary}

\begin{pf} $|\BL(\LL_Z)|\leq |\BL(\SS)|\leq 2^{2^{\aleph_0}}$, where the second inequality was proved in~\cite{[ohkawa]}.
\end{pf}

\begin{lemma}\label{compose}  Let $\T$ be as in Definition~\ref{cat}, and $X$ and $Y$ objects of $\T$.  Then $\<X\>\leq \<Y\>$ if and only if $L_X=L_XL_Y=L_YL_X$ and in this case the following diagram commutes (also with $\BL$ replaced by $\DL$ or $\BA$).

$$
\UseComputerModernTips
\xymatrix{ 
 & \BL(\T) \ar@{->>}[dl]_{L_Y} \ar@{->>}[dr]^{L_X}  \\
\BL(\LL_Y) \ar@{->>}[rr]^{L_X} && \BL(\LL_X)
}
$$
\end{lemma}

\begin{pf} The first equivalence is straightforward; it follows from~\cite[Prop.~1.22]{[Rav84]} and the observation that $\<X\>\leq \<Y\>$ precisely when all $L_X$-locals are $L_Y$-locals.  The last remark follows from Proposition~\ref{lattice-map}. \end{pf}

Here we mention one object-level application of the Bousfield lattice calculations of Sections~\ref{sect-enkn} and \ref{sect-other}.  Call an object $X$ {\it square-zero} if it is nonzero, but $X\w X=0$.  For example, $I\w I=0$ in $\SS$.

\begin{proposition} There are no square-zero objects in the harmonic, $E(n)$-, $K(n)$-, $\hfp$-, or $I$-local categories. \label{no-sz}
\end{proposition}

\begin{pf} In Corollary 2.8 of~\cite{[preprint]}, we show that in a well generated tensor triangulated category, there are no square-zero objects if and only if $\BA=\DL=\BL$.  The claim then follows from Corollaries~\ref{ln-bl} and~\ref{bl-kn} and Propositions~\ref{bl-hfp} and~\ref{bl-i}.
\end{pf}

\subsection{Background on spectra}  We quickly review some relevant background on the stable homotopy category.  See~\cites{[Rav92],[Hov95a],[MS], [RavLaTC]} for more details.  Fix a prime $p$ and let $\SS$ denote the $p$-local stable homotopy category of spectra.  Let $S^0$ denote the sphere spectrum.  The {\it finite} spectra $\FF$ are the compact objects of $\SS$, and $\FF = \th(S^0)$.  The structure of $\FF$ is determined by the Morava $K$-theories $K(i)$.  For each $i\geq 0$, $K(i)$ is a skew field object in $\SS$, such that $K(i)\w K(j)=0$ when $i\neq j$.  If $X$ is a finite spectrum and $K(j)\w X=0$, then $K(j-1)\w X=0$.  We say a finite spectrum $X$ is {\it type $n$} if $n$ is the smallest integer such that $K(n)\w X\neq 0$.  Define $\CC_n = \<K(n-1)\>\cap \FF$.  Then every thick subcategory of $\FF$ is $\CC_n$ for some $n$.  It follows that any two spectra of type $n$ generate the same thick subcategory, and hence Bousfield class; let $F(n)$ denote a generic type $n$ spectrum.

Given a type $n$ spectrum $X$, there is a $v_n$ self-map $f:\Sigma^d X\ra X$ with $[S^0, K(n)\w f]_i$ an isomorphism for all $i$ and $[S^0, K(m)\w f]_j=0$ for all $j$ and $m\neq n$.  We define $f^{-1}X$ to be the telescope, i.e.\ sequential colimit, i.e.\ homotopy colimit, of the diagram  $X\ra \Sigma^{-d} X\ra \cdots$.  By the periodicity theorem, any choice of $v_n$ self-map $f$ yields an isomorphic telescope.  The telescopes of different type $n$ spectra are Bousfield equivalent; denote this class by $\<T(n)\>$.

As mentioned above, localization away from a finite spectrum $F(n+1)$ exists, and is smashing.  This localization functor is denoted $\Lnf$, and is the same as homological localization at $T(0)\vee\cdots\vee T(n)$.

Let $E(n)$ denote the Johnson-Wilson spectrum; this is a ring spectrum with $\<E(n)\> = \<K(0)\vee\cdots\vee K(n)\>$.  Define $\Ln:\SS\ra \SS$ to be homological localization at $E(n)$.  A deep theorem of Ravenel~\cite[Thm.~7.5.6]{[Rav92]} shows that $\Ln$ is smashing for all $n$.  The $\Lnf$ and $\Ln$ are the only known smashing localization functors on $\SS$.  

The $\Lnf$-acyclics are $\loc(F(n+1)) = \loc(\CC_{n+1}) = \loc(\<K(n)\>\cap \FF) = \loc(\<E(n)\>\cap \FF)$.  The $\Ln$-acyclics are $\<E(n)\> = \loc(\<E(n)\>)$.  Thus every $\Lnf$-acyclic is $\Ln$-acyclic, and we have $\<K(0)\vee\cdots\vee K(n)\>\leq \<T(0)\vee\cdots \vee T(n)\>$ for all $n$.  It follows that there is a natural map $\Lnf \ra \Ln$.

For convenience later, we collect some calculations in $\SS$.

\begin{lemma} \label{calcs} In $\BL(\SS)$ we have the following.
\begin{enumerate}
    \item $\<F(m)\>\leq \<F(n)\>$ if and only if $m\geq n$.  For all $n$ and $m$, $\<F(m)\w F(n)\>\neq \0$.  Furthermore, $\<F(n)\w F(n)\> = \<F(n)\>$ for all $n$.
		\item $\<F(m) \w T(n)\> = \0$ when $m> n$, and $\<F(m) \w T(n)\> = \<T(n)\>$ when $m\leq n$.
		\item $\<T(m) \w T(n)\> = \0$ when $m\neq n$, and $\<T(n) \w T(n)\> = \<T(n)\>$.
		\item $\<F(m) \w K(n)\> = \0$ when $m > n$, and $\<K(n)\> = \<F(m) \w K(n)\> \leq \<F(m)\>$ when $m\leq n$.
		\item $\<T(m) \w K(n)\> = \0$ when $m \neq n$, and $\<K(n)\> = \<T(n) \w K(n)\> \leq \<T(n)\>$.
		\item $\<K(m) \w K(n)\> = \0$ when $m \neq n$, and $\<K(n) \w K(n)\> = \<K(n)\>$.
	\end{enumerate}
\end{lemma}

\begin{pf} Part $(1)$ is Theorem $14$ of~\cite{[HS]}, along with the observation that $\<F(n)\>$ is complemented by $\<L_{n-1}^fS^0\>$ and hence is in $\DL$.

Part (2) is in ~\cite[2.8(i)]{[RavLaTC]}, \cite[6.2]{[MS]}, and \cite[Sect.~5]{[HP]}.  Part $(3)$ is also in~\cite[Sect.~5]{[HP]}.

Part (4) follows from the definition of type $m$ spectra.  Since each $K(i)$ is a skew field object, $F(m)\w K(i)\neq \0$ implies this $K(i)$-module object $F(m)\w K(i)$ is a wedge of suspensions of $K(i)$.

From the periodicity theorem, $T(m)$ has nonzero $K(m)$ homology, so $T(m)\w K(m)\neq 0$.  The rest of Part $(5)$ is in~\cite[Prop.~A.2.13]{[Rav92]}.  Finally, Part $(6)$ is well known.
\end{pf}

\vspace{22pt}
\section{Local versions of the telescope conjecture} \label{sect-TCs}
\vspace{22pt}

In this section, let $L=L_Z:\SS\ra \SS$ be a localization functor for some $Z\in\SS$, and let $\LL=\LL_Z$ denote the category of $L$-locals.  First we state the various versions of the original telescope conjecture on $\SS$.

\begin{definition} Fix an integer $n\geq 0$. On $\SS$, we have the following versions of the telescope conjecture.\\

$\begin{array}{ll}
\TCi: & \<T(n)\> = \<K(n)\>. \\
\TCii: & L_n^f X\stackrel{\sim}{\ra} L_n X \mbox{ for all } X. \\
\TCiii: & \mbox{If }X\mbox{  is type }n\mbox{  and }f\mbox{  is a }v_n\mbox{  self-map, then }L_nX\cong f^{-1}X.\\
\GSC: & \mbox{Every smashing localization is generated by a set of compact objects.}\\
\SDGSC: & \mbox{Every smashing localization is generated by a set of }\\
 & \mbox{     strongly dualizable objects.}\\
\end{array}$
\end{definition}

\begin{theorem} \label{TCs} On the category $\SS$ we have $\TCi \iff \TCiii.$

Also, $\TCii$ holds if and only if $\TCij$ holds for all $i\leq n$.

Given $\TCiijj$ and $\TCi$, then $\TCii$ holds.

Furthermore, $\GSC \iff \SDGSC$, and this implies $\TCii$ for all $n$.
\end{theorem}

\begin{remark} Note that if we quantify over all $n$, the first three versions of the telescope conjecture are equivalent.  That is,
$$ \TCi \mbox{ for all }n \iff\TCii \mbox{ for all }n \iff \TCiii \mbox{ for all }n .$$
\end{remark}

\begin{remark} \label{just-s} $\TCii$ holds if and only if $L_n^f S^0\stackrel{\sim}{\ra} L_n S^0$.  Indeed, since both $\Lnf$ and $\Ln$ are smashing, the subcategory of objects $W$ such that $\Lnf W \stackrel{\sim}{\ra} \Ln W$ is localizing.  Thus if it contains $S^0$, it contains $\loc(S^0) = \SS$.
\end{remark}

\begin{pf} First we show the equivalence of $\TCi$ and $\TCiii$.  This is also sketched in~\cite[1.13]{[MRS]}.  For any type $n$ spectrum $Y$, $\th(Y) = \th(F(n))$ and so we have $\th(\Ln Y) = \th (\Ln F(n))$, and $\<\Ln Y\> = \<\Ln F(n)\>$.  A construction in~\cite[8.3]{[Rav92]} gives a type $n$ spectrum $Y$ with $\Ln Y\in \th(K(n))$.  Thus $\<\Ln Y\>\leq \<K(n)\>$, and $0\neq \<\Ln F(n)\> = \<K(n)\>$.  Suppose $\TCiii$ holds.  Then $\<\Ln Y\> = \<f^{-1} Y\> = \<T(n)\>$, and so $\<T(n)\> = \<K(n)\>$.

If $X$ is type $n$ and $f$ is a $v_n$ self-map, then~\cite[Prop.~3.2]{[MS]} implies that $\Lnf X\cong L_{T(n)}X\cong f^{-1}X$.  Thus assuming $\TCi$, we have $L_{K(n)}X\cong f^{-1}X$, and to prove $\TCiii$ it suffices to show that $\Ln X\cong L_{K(n)}X$.  This is known (see e.g.~\cite{[Hov95a]}), but we will give a proof that extends well to the localized setting.  Since $\<K(n)\>\leq \<E(n)\>$, localization at $K(n)$ gives a map $\Ln X\ra L_{K(n)}X$.  It suffices to show that this is an $\Ln$-equivalence.  The fiber is $K(n)$ acyclic, so  $\Ln X\w K(n) \ra L_{K(n)}X\w K(n)$ is an isomorphism.  Consider $i<n$.  The triangle $C_n X \w K(i)\ra X\w K(i)\ra \Ln X\w K(i)$ shows that $\Ln X\w K(i)$ is zero, because $X$ is type $n$ and $C_n X$ is $K(i)$ acyclic.  Lemma $3.3.1$ in~\cite{[HPS]} states that $LW = LS^0\w W$ for any localization $L$ and strongly dualizable $W$.  Since every finite spectrum is strongly dualizable, $L_{K(n)}X\w K(i) = L_{K(n)}S^0\w X\w K(i) = 0$.  Thus $\Ln X\w K(i) \ra L_{K(n)}X\w K(i)$ is an isomorphism for all $i\leq n$, and hence $\Ln X\ra L_{K(n)}X$ is an $\Ln$-equivalence.

For the second statement, note that $\TCii$ is equivalent to the statement $\<T(0)\vee\cdots\vee T(n)\> = \<K(0)\vee\cdots\vee K(n)\>$.  Smashing this with $\<T(i)\>$, for $0\leq i\leq n$, and using Lemma~\ref{calcs}, yields $\TCij$ for each $i$.  The third statement is also clear from this observation.

Finally, $\GSC\iff \SDGSC$ because objects in $\SS$ are compact if and only if they are strongly dualizable.  Given $\GSC$, consider $\Ln$.  The $\GSC$ would imply that the $\Ln$-acyclics are $\loc(\<E(n)\>\cap\FF)$.  As observed earlier, this is the same as $\loc(F(n+1))$, which are the $\Lnf$-acyclics.  Therefore we would have $\Lnf\cong \Ln$. \end{pf}

The $\SDGSC$ is new, and we will discuss it first.  As mentioned in the above proof, in $\SS$ compactness is equivalent to strong dualizability.  It is well known that localization away from a set of compact objects is smashing.  The $\GSC$ is precisely the statement that the converse holds.  However, as we will show next, one only needs strong dualizability to generate a smashing localization functor.  We will state our result in slightly more general terms.

\begin{theorem} \label{sdsl} Let $\T$ be a well generated tensor triangulated category such that  $\loc(\one) = \T$, as in Definition~\ref{cat}.  Let $A=\{B_\alpha\}$ be a (possibly infinite) set of strongly dualizable objects.  Then there exists a smashing localization functor $L:\T\ra \T$ with Ker~$L = \loc(A)$.
\end{theorem}

\begin{pf}  Let $E=\vee_\alpha B_\alpha$ and note that $\loc(E) = \loc(A)$.  The category $\T$ is well generated by hypothesis.  The localizing subcategory $\S=\loc(E)$ is also well generated, by~\cite[Rmk.~2.2]{[IK]}, and is tensor-closed by Lemma~\ref{useful}.

By~\cite[Prop.~2.1]{[IK]} there exists a localization functor $L:\T\ra \T$ with Ker~$L= \S$.  We will show that $L$ is a smashing localization.

First we claim that the $L$-locals are tensor-closed.  For any $Y\in \T$, we have the following.
$$\begin{array}{ll}
Y \mbox{ is }L\mbox{-local} & \iff [W,Y]_n=0 \mbox{ for all }W\in\S \mbox{ and all }n\in\bz\\
 & \iff [E, Y]_n=\prod[B_\alpha, Y]_n=0 \mbox{ for all }n\in\bz\\
 & \iff [B_\alpha, Y]_n=0 \mbox{ for all }\alpha \mbox{ and }n\in\bz\\
 & \iff DB_\alpha\w Y=0 \mbox{ for all }\alpha\\
\end{array}$$

The second equivalence follows from the fact that $\{X\;|\; [X,Y]_n=0 \mbox{ for all }n\in\bz\}$ is a localizing subcategory containing $E$, hence all of $\S$.  The final equivalence uses the fact that the $B_\alpha$ are strongly dualizable.

Now suppose $Y$ is $L$-local and $X$ is arbitrary.  Then $DB_\alpha\w Y=0$ for all $\alpha$, so $DB_\alpha\w Y\w X=0$ for all $\alpha$, and thus $Y\w X$ is $L$-local.  This shows that the $L$-locals are tensor-closed.

Consider the localization triangle $C\one\ra \one \ra L\one$, where $L\one$ is $L$-local and $C\one\in\S$.  For arbitrary $X\in \T$, tensoring gives an exact triangle
$$ C\one \w X \ra X\ra L\one\w X.$$

The object $L\one\w X$ is $L$-local, since the locals are tensor-closed.  Likewise, $C\one\w X\in \S$, since $\S$ is tensor-closed and so $L(C\one\w X)=0$.  Therefore $X\ra L\one \w X$ is an $L$-equivalence from $X$ to an $L$-local object, and it follows that $LX\cong L\one\w X$.  This shows that $L$ is a smashing localization. \end{pf}

In the  stable homotopy category, and more generally whenever $\one\in \T$ is compact, this gives nothing new; by~\cite[Thm.~2.1.3(d)]{[HPS]} compact and strongly dualizable are equivalent.  Consider, however, the harmonic category, which has no nonzero compact objects.  In Section~\ref{sect-harm} we classify all smashing localizations on the harmonic category; they are indexed by $\bn$.  Thus the $\GSC$ fails in the harmonic category, but as we show in Theorem~\ref{harm-sdgsc}, the $\SDGSC$ holds.  In fact, in the following sections we will give several examples of categories where the $\GSC$ fails but the $\SDGSC$ holds. 

On the other hand, we don't expect the $\SDGSC$ to hold in complete generality, since Keller's counterexample~\cite{[keller]} to the $\GSC$ is also a counterexample to the $\SDGSC$; in the derived category of a ring $R$, the unit $R$ is compact and strongly dualizable, so the $\GSC$ and $\SDGSC$ are equivalent.

The $\GSC$ and $\SDGSC$ make sense in any localized category, but $\TCi$, $\TCii$, and $\TCiii$ may not, since $T(n)$ and $K(n)$ may not be objects in $\L$.  Instead we make the following definitions.

\begin{definition} Let $L:\SS\ra \SS$ be a localization, and $\LL$ the category of $L$-locals. \label{defn-l}
 \begin{enumerate}  
	\item Let $\lnf:\LL\ra \LL$ denote localization at $\<LT(0)\vee LT(1)\vee\cdots \vee LT(n)\>$.
	 \item Let $\ln: \LL\ra \LL$ denote localization at $\<LK(0)\vee LK(1)\vee\cdots \vee LK(n)\>$.
	\end{enumerate}
\end{definition}

Before stating and proving a local version of Theorem~\ref{TCs}, we establish some results about $\lnf$ and $\ln$.  First we make an observation about calculations in $\BL(\LL)$.%In order to do this, it is useful to collect some calculations.

\begin{lemma} \label{L-calcs} All the calculations in Lemma~\ref{calcs} are valid in $\BL(\LL)$ if we replace $F(n)$, $T(n)$, and $K(n)$ with $LF(n)$, $LT(n)$, and $LK(n)$.
\end{lemma}

\begin{pf} This follows from Theorem~\ref{local-cat} and the statements in Lemma~\ref{calcs}.
\end{pf}

\begin{proposition} \label{lnf-smash} $\lnf$ is localization away from $LF(n+1)$, and hence is smashing.
\end{proposition}

\begin{pf} By Theorem~\ref{sdsl} we know that there is some smashing localization functor $\l: \LL\ra \LL$ that is localization away from $LF(n+1)$; we wish to show $\l=\lnf$.  Let $\one=LS^0$ for simplicity of notation, and let $\c$ denote the colocalization corresponding to $\l$.  We claim that the $\l$-acyclics are precisely $\loc(\c\one)$.  Clearly $\c\one$ is $\l$-acyclic, and these are a localizing subcategory, so $\loc(\c\one)\subseteq \{\l\mbox{-acyclics}\}$.  On the other hand, suppose $W$ is $\l$-acyclic.  Because $\l$ is smashing, $LW = W\w \l\one=0$, so $W=W\w c\one$.  Then since $W\in \loc(\one) = \LL$, we have $W=W\w \c\one\in\loc(\one\w\c\one)=\loc(\c\one)$, proving the claim.

By definition, the $\l$-acyclics are also given by $\loc(LF(n+1))$.  Therefore $\<LF(n+1)\> = \<\c\one\>$.  

The class $\<F(n+1)\>$ is complemented by $\<T(0)\vee\cdots\vee T(n)\>$ in $\BL(\SS)$~\cite[Sect.~5]{[HP]}, so $\<LF(n+1)\>$ is complemented by $\<LT(0)\vee\cdots\vee LT(n)\>$ in $\BL(\LL)$.  At the same time, $\<\c\one\>$ is complemented by $\<\l\one\>$, and complements are unique.  We conclude that $\{\l\mbox{-acylics}\} = \<\l\one\> = \<LT(0)\vee\cdots\vee LT(n)\>$.  Since $\l$ and $\lnf$ are two localizations on $\LL$ with the same acyclics, they are equal. \end{pf}

\begin{lemma} \label{sm-compose} If $L$ is smashing, then $\lnf = L\Lnf=\Lnf L$ and $\ln = L\Ln=\Ln L$, and both are smashing.
\end{lemma}

\begin{pf} Smashing localization functors always commute, and compose to give a smashing localization.  The functor $L\Ln:\SS\ra\SS$, sending $X\mapsto L(\Ln S^0\w X) =  LS^0\w \Ln S^0 \w X$ is a smashing localization.  Since $L\Ln$-locals are $L$-local, it also gives a smashing localization on $\LL$.  The acyclics of this functor are $\<L\Ln S^0\>$ in $\BL(\LL)$, which is $\<LK(0)\vee\cdots\vee LK(n)\>$.  Thus $L\Ln$ and $\ln$ are localizations on $\LL$ with the same acyclics, hence isomorphic.  The same proof works for $\lnf = L\Lnf$. \end{pf}

In $\SS$, for a type $n$ finite spectrum $X$ with a $v_n$ map $f:\Sigma^dX\ra X$ and telescope $f^{-1}X$, it is known~\cite[Prop.~3.2]{[MS]} that $\Lnf X\cong L_{T(n)}X\cong f^{-1}X$.  The following proposition shows that the local version of this result holds as well.

\begin{lemma} \label{ugly} Let $L:\SS\ra \SS$ be a localization, and $\lnf$, $X$ and $f^{-1}X$ as above.  Then $$\lnf (LX)\cong L_{LT(n)}(LX)\cong L(f^{-1}X).$$
\end{lemma}

\begin{pf}  The proof parallels the~\cite{[MS]} result, one must only check that everything works when localized.  If $L_{LT(n)}(LX)\cong L(f^{-1}X)$ holds for a single type $n$ spectrum, then it holds for all type $n$ spectra.  So without loss of generality, we can choose $X$ to be a type $n$ spectrum that is a ring object in $\SS$.  Then for any $v_n$ self-map $f$, the telescope $f^{-1}X$ is also a ring object~\cite[Lemma~2.2]{[MS]}.  By Lemma~\ref{useful}, $L(f^{-1}X)$ is a ring object in $\LL$, and hence is local with respect to itself.  

Lemma~2.2 in~\cite{[MS]} shows that $X\w f^{-1}X \cong f^{-1}X\w f^{-1}X$ in $\SS$, so $LX \w_\LL L(f^{-1}X) \cong L(f^{-1}X)\w_\LL L(f^{-1}X)$ in $\LL$ and the canonical map $LX\ra L(f^{-1}X)$ is an $L(f^{-1}X)$-equivalence.  It follows that $L_{LT(n)}(LX)\cong L(f^{-1}X)$.

Since $\<LT(n)\>\leq \<LT(0)\vee\cdots\vee LT(n)\>$, $L(f^{-1}X)$ is $\lnf$-local.  One then uses Lemma~\ref{L-calcs} to see that $LX\ra L(f^{-1}X)$ is a $\lnf$-equivalence, and so $\lnf(LX) = \lnf(L(f^{-1}X))=L(f^{-1}X)$. \end{pf}

\begin{definition} Let $L:\SS\ra \SS$ be a localization, and consider the category $\L$ of locals.  Fix an $n\geq 0$.  We have the following versions of the telescope conjecture on $\L$.  \\

$\begin{array}{ll}
\LTCi: & \<LT(n)\> = \<LK(n)\>. \\
\LTCii: & \lnf X\stackrel{\sim}{\ra} \ln X \mbox{ for all } X. \\
\LTCiii: & \mbox{If }X\in \SS \mbox{  is type }n\mbox{  and }f\mbox{  is a }v_n\mbox{  self-map, then }\ln (LX)\cong L(f^{-1}X).\\
\GSC: & \mbox{Every smashing localization is generated by a set of compact objects.}\\
\SDGSC: & \mbox{Every smashing localization is generated by a set of }\\
 & \mbox{     strongly dualizable objects.}\\
\end{array}$
\end{definition}

\begin{theorem} \label{LTCs} On the category $\L$ we have $\LTCi\Rightarrow \LTCiii.$

Also, $\LTCii$ holds if and only if $\LTCij$ holds for all $i\leq n$.

Given $\LTCiijj$ and $\LTCi$, then $\LTCii$ holds.
\end{theorem}

\begin{pf} Note that $\LTCii$ is equivalent to the statement $\<LT(0)\vee\cdots\vee LT(n)\> = \<LK(0)\vee\cdots\vee LK(n)\>$, so the last two statements are clear.  We will show that $\LTCi\Rightarrow \LTCiii$, by mimicking the proof in Theorem~\ref{TCs}.

If $X$ is type $n$ and $f$ is a $v_n$ self-map, Lemma~\ref{ugly} shows that $\lnf (LX)\cong L_{LT(n)}(LX)\cong L(f^{-1}X)$.  Then $\LTCi$ implies $L_{LK(n)}LX\cong L(f^{-1}X)$.  So it suffices to show that $\ln (LX) = L_{LK(n)} LX$.  We must show that the map $L_{LK(n)}: \ln (LX) \ra L_{LK(n)} LX$ is an $\ln$-equivalence.  The same reasoning as in Theorem~\ref{TCs}, along with the computations of Lemma~\ref{L-calcs} and some definition unwinding, gives us that $\ln (LX)\w LK(i) \ra L_{LK(n)} LX\w LK(i)$ is an equivalence for all $i\leq n$; we only need to notice that Lemma $3.3.1$ in~\cite{[HPS]} applies to strongly dualizable objects, and $LX$ is strongly dualizable.  \end{pf}

\begin{theorem} \label{LTCstwo} Furthermore, if $L:\SS\ra \SS$ is a smashing localization, then on the category $\L$ of locals we have $\LTCi\Leftarrow \LTCiii,$ and
$$\GSC \iff \SDGSC \Rightarrow \LTCii \mbox{ for all }n.$$
\end{theorem}

\begin{remark} In this case, $\LTCii$ is equivalent to $\lnf (LS^0)\stackrel{\sim}{\ra} \ln (LS^0)$, since by Lemma~\ref{sm-compose} both $\lnf$ and $\ln$ are smashing, so the argument in Remark~\ref{just-s} applies.
\end{remark}

\begin{pf}  By~\cite[Thm.~2.1.3(d)]{[HPS]}, the compact objects and strongly dualizable objects in $\LL$ coincide.  Thus $\GSC\iff\SDGSC$, and this implies $\LTCii$ just as in Theorem~\ref{TCs}.

Suppose $X$ has type $n$ and $\LTCiii$ holds.  As in the proof of Theorem~\ref{TCs}, $\<\Ln F(n)\>=\<K(n)\>$ in $\BL(\SS)$, so $\<L\Ln LF(n)\> = \<LK(n)\>$ in $\BL(\LL)$.  By Lemma~\ref{sm-compose}, we have $\<\ln LF(n)\> = \<LK(n)\>$.  Then $\LTCiii$ implies that $\<LT(n)\> = \<L(f^{-1}X\> = \<\ln(LX)\> = \<\ln LF(n)\>$, so $\LTCi$ holds.
\end{pf}

\begin{question} Is $\ln:\LL\ra\LL$ always a smashing localization?
\end{question}

This is the case in all the local categories investigated in this paper, whether or not $L:\SS\ra \SS$ is a smashing localization.  If one could show $\ln$ is always smashing, then most likely on $\LL$ one would have $\SDGSC\implies \LTCii \forall n$.

We would of course like to know if and when information on localized telescope conjectures can help with those in the original category $\SS$, where all versions remain open.  

\begin{proposition} \label{implications} Let $L:\SS\ra \SS$ be a localization, with localized category $\LL$.  
\begin{enumerate}
  \item If $\TCi$ holds on $\SS$, then $\LTCi$ holds on $\LL$.
  \item If $\TCii$ holds on $\SS$, then $\LTCii$ holds on $\LL$.
	\item If $\TCiii$ holds on $\SS$, then $\LTCiii$ holds on $\LL$.
\end{enumerate}

\noindent Furthermore, if $L$ is a smashing localization then we have the following.

\begin{enumerate}
	\item[(4)] If $\GSC$ holds on $\SS$, then $\GSC$ holds on $\LL$.
	\item[(5)] If $\SDGSC$ holds on $\SS$, then $\SDGSC$ holds on $\LL$.
\end{enumerate}
\end{proposition}

\begin{pf} Part $(1)$ follows immediately from Proposition~\ref{lattice-map}.  So does Part $(2)$, since $\TCii$ is equivalent to the statement $\<T(0)\vee\cdots\vee T(n)\>=\<K(0)\vee\cdots\vee K(n)\>$ and similarly for $\LTCii$.  From this and Theorems~\ref{TCs} and~\ref{LTCs} we have $\TCiii\iff \TCi\Rightarrow \LTCi\Rightarrow  \LTCiii.$

Now suppose $L$ is smashing, and the $\GSC$ holds on $\SS$.  Let $\l:\LL\ra \LL$ be a smashing localization.  Thus $\l$ is defined by $\l(LY) = \l(LS^0)\w_\LL LY = \l S^0\w LS^0\w Y$.  We can therefore extend $\l$ to be a smashing localization on all of $\SS$, with $X\mapsto \l S^0\w LS^0\w X = \l LS^0\w X$.  Since the $\GSC$ holds on $\SS$ by assumption, the acyclics of this functor are $\<\l LS^0\> = \loc(A)$, for some set of compact objects $A$ in $\SS$.  Here $\<\l LS^0\>$ refers to the Bousfield class in $\BL(\SS)$.

We must show that $\<\l LS^0\>$ in $\BL(\LL)$ is generated by a set of objects that are compact in $\LL$.  Note that $\<\l LS^0\>$ in $\BL(\LL)$ is $\{LW\;|\; LW\w_\LL \l LS^0=0\} = \{LW\;|\; LW\w_\SS \l LS^0=0\} = \<\l LS^0\> \cap \LL$, where the latter $\<\l LS^0\>$ is in $\BL(\SS)$.  Therefore $\<\l LS^0\>$ in $\BL(\LL)$ is $\loc(A)\cap \LL$.  We claim that this is $\loc(L(A))$.  Since $L$ sends compacts to compacts, this will show that $\l$ is generated by a set of compacts.  

If $X\in \loc(A)$, then $LX\in \loc(L(A))$.  If $X\in \LL$ in addition, then $X\cong LX\in \loc(L(A))$.  For the other inclusion, note that the intersection of two localizing subcategories is a localizing subcategory, and $\LL$ is a localizing subcategory of $\SS$ because $L$ is smashing.  For $Y\in A$, $LY$ is in $\LL$, and $LY = LS^0\w Y\in \loc(Y)\subseteq \loc(A)$.  Therefore $L(A)\subseteq \loc(A)\cap \LL$, and $\loc(L(A))=\loc(A)\cap \LL$.

Part $(5)$ follows immediately, since if $L$ is smashing then $\GSC\iff \SDGSC$ in both $\SS$ and $\LL$. \end{pf}

Balmer and Favi~\cite[Prop.~4.4]{[BF]} have also recently proved Part $(4)$ in the slightly more general setting of a smashing localization on a unital algebraic stable homotopy category; the above proof would apply there as well.  One would like to prove Part $(5)$ without the assumption that $L$ is smashing, but it's not clear if this is possible.

Letting $L=L_Z$, for $Z = \bigvee_{i\geq 0}K(i)$, $E(n)$, $K(n)$, $BP$, $H\bf_p$, or $I$ provides interesting examples of categories $\L$ on which to investigate these telescope conjectures.  Furthermore, $\LTCi$ suggests the relevance of Bousfield lattices to understanding these questions.  In the remaining sections, we focus on specific localized categories.

\vspace{22pt}
\section{The Harmonic category} \label{sect-harm}
\vspace{22pt}

Let $Q=\bigvee_{i\geq 0}K(i)$, and $L=L_Q:\SS\ra \SS$, and consider the harmonic category $\HH$ of $L$-locals.  Harmonic localization is not smashing.  An object is called {\it harmonic} if it is $L$-local, and {\it dissonant} if it is $L$-acyclic.  For example, finite spectra, suspension spectra, finite torsion spectra, and $BP$ are known to be harmonic~\cites{[Hov95a],[Rav84]}.  On the other hand, $I$ and $\hfp$ are dissonant.

In order to answer the telescope conjectures in $\HH$, we will first calculate the Bousfield lattice of $\HH$.  In this section all smash products are in $\HH$, unless otherwise noted.  Given any set $P$, let $2^P$ denote the power set of $P$.

\begin{definition} Given $X\in \HH$, define the {\it support} of $X$ to be $$\supp(X) = \{i\;|\; X\w K(i)\neq 0\}\subseteq \bn.$$
\end{definition}

The following result and proof was pointed out to us by Jon Beardsley.

\begin{proposition} \label{jb} The Bousfield lattice of $\HH$ is $2^\bn$.
\end{proposition}

\begin{pf} Each $K(n)$ is a ring object, hence $K(n)$-local by Lemma~\ref{useful}.  Because $\<K(n)\>\leq \<Q\>$, $K(n)$-locals are harmonic, thus each $K(n)$ is harmonic.  The argument hinges on the fact that $K(n)$ is a skew field object in $\HH$: for $X=LX$ in $\HH$, if $X\w K(n)\neq 0$ then $X\w K(n) = L(X\w_\SS K(n))$ so $X\w_\SS K(n)\neq 0$, and $X\w_\SS K(n)$ must be a nonempty wedge of suspensions of $K(n)$'s.  Thus
$$X\w K(n) = L(X\w_\SS K(n)) = L(\vee \Sigma^i K(n)) = L(\vee \Sigma^i LK(n)) = \coprod_\LL \Sigma^i LK(n) = \coprod_\LL \Sigma^i K(n).$$

It follows that $LX\w K(n)=0$ if and only if $LX\w_\SS K(n)=0$.  Furthermore, if $LX\w K(n)\neq 0$, then $\<LX\w K(n)\> = \<K(n)\>$, where these are Bousfield classes in $\BL(\HH)$.

By the definition of $L$, for any $W\in \SS$, if $W\w_\SS K(n) =0$ for all $n$, then $LW=0$.  Combining this with the above observation, we get that a local object $W=LW$ has $W\w K(n)=0$ in $\HH$ for all $n$ if and only if $W=0$.

Therefore, for any $X,Y\in\HH$, we have 
$$Y\w X=0 \iff Y\w X\w K(n)=0 \;\forall n \iff Y\w K(n)=0 \mbox{ for all }n\in \supp(X).$$

We conclude that there is a lattice isomorphism $F:\BL(\HH)\ra 2^\bn$, given by the following.
$$ \<X\> = \bigvee_{\supp(X)} \<K(i)\> \mapsto \supp(X),$$
$$ N\subseteq\bn \mapsto  \bigvee_{i\in N} \<K(i)\>.$$ \end{pf}

\begin{theorem} On $\HH$, for all $n\geq 0$ we have that $\LTCi, \LTCii$, and $\LTCiii$ hold. \label{LTCs-harm}
\end{theorem}

\begin{pf} By Lemmas~\ref{calcs} and~\ref{L-calcs}, $LT(n)$ and $LK(n) = K(n)$ have the same support.  The above theorem then implies that $\<LT(n)\> = \<LK(n)\>$.  Thus $\LTCi$ holds for all $n$, and the claim follows from Theorem~\ref{LTCs}.
\end{pf}

Next, we classify all smashing localizations on $\HH$, and show that the $\GSC$ fails but the $\SDGSC$ holds.  The proof is based on that of~\cite[Thm.~6.14]{[HovS]}, which classifies smashing localizations in the $E(n)$-local category.

\begin{theorem} \label{harm-sdgsc} If $L':\HH\ra \HH$ is a smashing localization functor, then $L'=\lnf$ for some $n\geq 0$, or $L'=0$ or $L'=id$.  Therefore the $\GSC$ fails but the $\SDGSC$ holds on $\HH$.
\end{theorem}

\begin{pf} Let $L':\HH\ra \HH$ be a smashing localization functor, and let $\one=LS^0$ be the unit in $\HH$.  The acyclics of $L'$ are given by $\<L'\one\>$.  From Proposition~\ref{jb}, $\<L'\one\>$ is equal to the wedge of $\<K(i)\>$ for all $i\in \supp(L'\one)$.  If $\supp(L'\one)=\emptyset$ then $\<L'\one\>=\0$ and $L'=0$.  

Assume now that $\supp(L'\one)$ is not empty, and take $j\in\supp(L'\one)$.  We will show that $\<L'\one\>\geq \<K(0)\vee\cdots\vee K(j)\>$.  It follows that either $\<L'\one\> = \bigvee_{i\geq 0}\<K(i)\> = \<\one\>$ and $L'=id$, or $L'=\ln=\lnf$ for $n=max(\supp(L'\one))$.

Since $\<K(j)\>\leq \<L'\one\>$, from Lemma~\ref{compose} we have $L_{K(j)}L' = L'L_{K(j)}=L_{K(j)}$.  Therefore $\<L_{K(j)}\one\> = \<L'\one\w L_{K(j)}\one\>\leq \<L'\one\>$.  Proposition 5.3 of~\cite{[HovS]} shows that in $\SS$, $L_{K(j)}S^0\w_\SS K(i)$ is nonzero for $0\leq i\leq j$ and zero for $i>j$.  Note that $L_{K(j)}S^0 = L_{K(j)}LS^0=L_{K(j)}\one$, and as remarked in the proof of Proposition~\ref{jb}, $LX\w K(i)=0$ if and only if $LX\w_\SS K(i)=0$.  Therefore in $\BL(\HH)$ we have $\<L_{K(j)}\one\> = \<K(0)\vee\cdots\vee K(j)\>$, and so $\<L'\one\>\geq \<K(0)\vee\cdots\vee K(j)\>$ as desired.

Each $\lnf$ is localization away from $LF(n+1)$ by Proposition~\ref{lnf-smash}, which is strongly dualizable by Theorem~\ref{local-cat}.  The identity is localization away from zero, and the zero functor is localization away from $LS^0$; these are both strongly dualizable.  Therefore the $\SDGSC$ holds.  On the other hand, Corollary B.13 in~\cite{[HovS]} shows that there are no nonzero compact objects in $\HH$, so the $\GSC$ fails. \end{pf}

\begin{question} Classify localizing subcategories of $\HH$.
\end{question}

It seems likely that every localizing subcategory of $\HH$ is a Bousfield class, and so these are in bijection with $2^\bn$, but we have been unable to prove this.

\vspace{22pt}
\section{The $E(n)$- and $K(n)$- local categories} \label{sect-enkn}
\vspace{22pt}

\subsection{The $E(n)$-local category} \label{sect-en}

Recall that $E(n)$ has $\<E(n)\> = \<K(0)\vee K(1)\vee\cdots \vee K(n)\>$.  In this section, fix $L=L_n = L_{E(n)}: \SS\ra \SS$ and let $\LLn$ denote the local category.  The functor $L_n$ is smashing, and so by Theorem~\ref{cat} each $LF(i)$ is compact in $\LLn$. Hovey and Strickland~\cite{[HovS]} have studied $\LLn$ in detail, and determine the localizing subcategories, smashing localizations, and Bousfield lattice of $\LLn$.  We begin by recalling these results.

\begin{lemma} For $0\leq i\leq n$, we have $LK(i) = K(i)$, and for $i >n$ we have $LK(i)=0$.
\end{lemma}

\begin{pf}  This follows from $\<E(n)\> = \<K(0)\vee K(1)\vee\cdots \vee K(n)\>$. \end{pf}

\begin{theorem}~\cite[Thm.~6.14]{[HovS]} The lattice of localizing subcategories of $\LLn$, ordered by inclusion, is in bijection with the lattice of subsets of the set $\{0,1,...,n\}$, where a localizing subcategory $\S$ corresponds to 
$$\{ i \;|\; K(i) \in \S\}.$$
\end{theorem}

\begin{corollary} Every localizing subcategory of $\LLn$ is a Bousfield class, in particular a localizing subcategory $\S$ is the Bousfield class $$\bigvee \< K(j)\;|\; K(j)\notin \S, 0\leq j\leq n\>.$$  
\end{corollary}

\begin{corollary} \label{ln-bl} For every $n\geq 0$, there is a lattice isomorphism $$f_n: \BL(\LLn)\stackrel{\sim}{\longrightarrow} 2^{\{0,1,...,n\}}.$$
\end{corollary}

\begin{pf} The isomorphism is given by
$$ \<X\> = \bigvee_{X\w K(i)\neq 0} \<K(i)\> \mapsto \{ i\;|\; X\w K(i)\neq 0\},$$
$$ N\subseteq \{0,1,...,n\} \mapsto  \bigvee_{i\in N} \<K(i)\>.$$
\end{pf}

\begin{theorem}~\cite[Cor.~6.10]{[HovS]} If $L':\LLn\ra \LLn$ is a smashing localization, then $L' = L_i=L_i^f$ for some $0\leq i\leq n$ or $L'=0$.  Thus the $\GSC$ holds on $\LLn$.
\end{theorem}

\begin{corollary} On $\LLn$, all of $\LTCij$, $\LTCiij$, $\LTCiiij$ hold for all $i$, and $\GSC$ and  $\SDGSC$ also hold. \label{LTCs-en}
\end{corollary}

\begin{pf} This follows from the previous Theorem, and Theorem~\ref{LTCstwo}.  Note that for $i>n$, we have $LT(i)=0=LK(i)$ by Lemma~\ref{calcs}, and so $\l_i=\ln = \lnf = \l_i^f$.
\end{pf}

Recall that there is a natural map $L_nX\ra L_{n-1}X$ for all $X$ in $\SS$ and $n$, and by Proposition~\ref{lattice-map} this induces a surjective lattice map $\BL(\LL_n)\ra \BL(\LL_{n-1})$, and an inverse system of lattice maps.

$$\cdots \ra \BL(\LL_n)\ra \BL(\LL_{n-1}) \ra \cdots \ra \BL(\LL_1)\ra \BL(\LL_{0})$$

\begin{proposition} The lattice isomorphisms $F$ and $f_n$ from Proposition~\ref{jb} and Corollary~\ref{ln-bl} realize $\BL(\HH)$ as the inverse limit of the maps $\BL(\LL_n)\ra \BL(\LL_{n-1})$. \label{realize}
\end{proposition}

\begin{pf} From Lemma~\ref{compose}, and the facts that $L_QK(i) = K(i)$ for all $i$, and $L_nK(i) = K(i)$ for $i\leq n$ and $L_nK(i)=0$ for $i>n$, we get the following diagram for all $n$.  The map $2^{\{0,1,...,n\}}\ra 2^{\{0,1,...,n-1\}}$ is induced by sending $m\mapsto m$ for $m<n$ but $n\mapsto 0$, and the maps $2^\bn \ra 2^{\{0,1,...,i\}}$ are defined similarly.

$$
\UseComputerModernTips
\xymatrix{ 
 & \BL(\HH) \ar@{->>}[dl] \ar@{->>}[dd] \ar@{<->}[rrr]_{F} &&& 2^\bn \ar@{->>}[dl] \ar@{->>}[dd] \\
\BL(\LL_n) \ar@{->>}[dr] \ar@{<->}[rrr]_{f_n} &&& 2^{\{0,1,...,n\}} \ar@{->>}[dr]\\
& \BL(\LL_{n-1}) \ar@{<->}[rrr]_{f_{n-1}} &&& 2^{\{0,1,...,n-1\}}\\
}
$$
\end{pf}

\subsection{The $K(n)$-local category}  Although an incredibly complicated category in its own right, the $K(n)$-local category is quite basic from the perspective of localizing subcategories, Bousfield lattices, and telescope conjectures.  In this subsection, let $L=L_{K(n)}: \SS\ra \SS$ be localization at $K(n)$, and let $\KKn$ denote the category of locals.  Hovey and Strickland classify the localizing subcategories of $\KKn$, and there are not many of them.\\

\begin{proposition}~\cite[Thm.~7.5]{[HovS]} There are no nonzero proper localizing subcategories of $\KKn$.
\end{proposition}

This Proposition implies that the Bousfield lattice of $\KKn$ is the two-element lattice $\{\<0\>, \<K(n)\>\}$.  We will prove a slightly more general result, that will be used again in Subsection~\ref{hfp}.

\begin{proposition} \label{sfobl} Consider a category $\T$ as in Definition~\ref{cat}, an object $Z$ in $\T$, and localization $L_Z:\T\ra \T$ with localized category $\LL_Z$.
\begin{enumerate}
   \item If $Z$ is a ring object, then $\<L_ZZ\>=\<Z\>$ is the maximum class in $\BL(\LL_Z)$.
	 \item If $Z$ is a skew field object, then $\BL(\LL_Z)$ is the two-element lattice $\{\0, \<Z\>\}$.
\end{enumerate}
\end{proposition}

\begin{pf} For $(1)$, first note that Lemma~\ref{useful} implies $L_ZZ=Z$.  Consider $\<Z\>$ in $\BL(\LL_Z)$.  By definition, this is the collection of all $W\in \LL_Z$ with $L(Z\w_\T W)=0$.  But $Z\w_\T W$ is a $Z$-module object in $\T$, hence is $L_Z$-local.  The only object that is both local and acyclic with respect to any localization is zero, so $Z\w_\T W=0$.  But this says that $W$ is $L_Z$-acyclic, hence zero in $\LL_Z$.  Therefore in $\BL(\LL_Z)$ we have $\<Z\> = \{0\}$.

Now suppose $Z$ is a skew field object in $\T$.  In particular, it is a ring object, so $\<Z\>$ is the maximum class in $\BL(\LL_Z)$.  Consider $\<LX\>$ in $\BL(\LL_Z)$, for arbitrary $X\in \T$.  If $X\w_\T Z=0$, then $LX=0$.  Otherwise, $X\w_\T Z$ is a wedge of suspensions of $Z$, so $\<Z\> = \<X\w_\T Z\>\leq \<X\>$ in $\BL(\T)$.  Then $\<Z\> = \<L_ZZ\>\leq \<L_ZX\>$ in $\BL(\LL_Z)$, so $\<L_ZX\>=\<Z\>$. \end{pf}

\begin{corollary} The Bousfield lattice of $\KKn$ is $\{\<0\>, \<K(n)\>\}$. \label{bl-kn}
\end{corollary}

\begin{theorem} In $\KKn$, all of $\LTCij$, $\LTCiij$, $\LTCiiij$ hold for all $i$, and $\GSC$ and  $\SDGSC$ also hold. \label{LTCs-kn}
\end{theorem}

\begin{pf} In light of Theorem~\ref{LTCs}, we will show that $\LTCij$ holds for all $i$.  This follows from Lemma~\ref{calcs}: for $i\neq n$ we have $LT(i)=0=LK(i)$, but $LT(n)\neq 0$ so by the last corollary $\<LT(n)\> = \<K(n)\> = \<LK(n)\>$.  

There are exactly two smashing localizations on $\KKn$.  The identity functor is smashing, and is localization away from $0$, which is compact and strongly dualizable.  The zero functor is smashing, and is localization away from $LS^0$, which is strongly dualizable.  It is not compact, but by Theorem 7.3 in~\cite{[HovS]}, $LF(n)$ is compact in $\KKn$ and $\loc(LF(n)) = \loc(LS^0)=\KKn$.  Therefore the zero functor is also generated by a compact object.  This shows that both the $\GSC$ and $\SDGSC$ hold. \end{pf}

\vspace{22pt}
\section{Other localized categories} \label{sect-other}
\vspace{22pt}

In this section we will consider several other localized categories.  In each case, let $L_Z:\SS\ra \SS$ denote the localization functor that annihilates $\<Z\>$, and let $\LL_Z$ denote the category of $L_Z$-locals.

\subsection{The $H\bf_p$-local category} \label{hfp} The Eilenberg-MacLane spectrum $\hfp$ is a skew field object in $\SS$; in fact, every skew field object in $\SS$ is isomorphic to either $\hfp$ or a $K(n)$.  Unlike the $\<K(n)\>$, it is not complemented; for example, $\<I\>\leq \<\hfp\>$ but $I\w \hfp=0$.  So $\<\hfp\>\in\DL\backslash\BA$.  Hovey and Palmieri~\cite{[HP]} have conjectured several results about the collection of classes less than $\<\hfp\>$ in $\BL(\SS)$.  The telescope conjectures and Bousfield lattice of $\LL_\hfp$ are quite simple.

\begin{theorem} In $\LL_{\hfp}$, all of $\LTCi$, $\LTCii$, $\LTCiii$ hold for all $n$. \label{LTCs-hfp}
\end{theorem}

\begin{pf} For all $n$, $K(n)\w \hfp = 0$ and $T(n)\w \hfp = 0$, by~\cite[p.~16]{[HP]}.  Therefore $LK(n)=0=LT(n)$ and $\LTCi$ holds for all $n$.  Note that $\ln = \lnf$ is the zero functor for all $n$.  \end{pf}

In order to discuss the $\GSC$ and $\SDGSC$ in this category, we must classify the smashing localizations.  We will do this by using what we know about the Bousfield lattice.

\begin{proposition} \label{bl-hfp} The Bousfield lattice of $\LL_\hfp$ is the two-element lattice $\{\0, \<\hfp\>\}$. 
\end{proposition}

\begin{pf} This follows immediately from Proposition~\ref{sfobl} because $\hfp$ is a skew field object in $\SS$.
\end{pf}

Recall that every smashing localization gives a pair of complemented classes in $\BA\subseteq \BL$.  Thus in $\LL_\hfp$ there are exactly two smashing localizations, the trivial ones given by smashing with zero and with the unit.

\begin{proposition} In $\LL_\hfp$, the $\GSC$ fails but the $\SDGSC$ holds. \label{gsc-hfp}
\end{proposition}

\begin{pf} The identity functor is smashing, and is localization away from $0$, which is compact and strongly dualizable.  By~\cite[Cor.~B.~13]{[HovS]}, there are no nonzero compact objects in $\LL_\hfp$.  So the zero functor, which is localization away from $LS^0$, is generated by a strongly dualizable object but not a compact one. \end{pf}

%\begin{question} Classify the localizing subcategories of $\LL_\hfp$.
%\end{question}
%
%The object $L\hfp = \hfp$ is also a skew field object in $\LL_\hfp$.  So $\loc(\hfp)$ is a minimal nonzero localizing subcategory in $\LL_\hfp$.  If one could show that $\loc(\hfp)\neq \LL_\hfp$, this would demonstrate a localizing subcategory that is not a Bousfield class, in a topological setting.  To date, the only known example of such a localizing subcategory is in the derived category of an absolutely flat ring that is not semi-artinian~\cite{[stevenson-counterexample]}.

One application of this Bousfield lattice calculation is to the question of classifying localizing subcategories.  Every Bousfield class is a localizing subcategory.  In Conjecture 9.1 of~\cite{[HP]}, Hovey and Palmieri conjecture the converse holds, in the $p$-local stable homotopy category.  The original conjecture is still open, but the question can be asked in any well-generated tensor triangulated category.  For example, in a stratified category every localizing subcategory is a Bousfield class.  The question is interesting, since in general localizing subcategories are hard to classify.  In many cases, including $\SS$, it is not even known if there is a set of localizing subcategories.  Recently Stevenson~\cite{[stevenson-counterexample]} found the first counterexample, in an algebraic setting: in the derived category of an absolutely flat ring that is not semi-artinian, there are localizing subcategories that are not Bousfield classes.  Now we show that $\LL_\hfp$ provides another counterexample.

\begin{proposition} \label{loc-not-bc} In $\LL_\hfp$ there are localizing subcategories that are not Bousfield classes.
\end{proposition}

\begin{pf} The following counterexample was suggested to us by Mark Hovey.  The Bousfield lattice of $\LL_\hfp$ has only two elements: $\0 = \LL_\hfp$ and $\one = \{0\}$.  It suffices to find a proper nonzero localizing subcategory in $\LL_\hfp$.  

Consider the Moore spectrum $M(p)$, defined by the triangle $S^0 \stackrel{p}{\ra} S^0 \ra M(p)$; this spectrum is $\hfp$-local.  Consider the following full subcategory in $\LL_\hfp$.

$$\A = \{ X\in \LL_\hfp\;|\; [X, M(p)]_n=0 \mbox{ for all } n\in\bz\}.$$

This is a localizing subcategory, called the cohomological Bousfield class of $M(p)$ and denoted $\<M(p)^*\>$ in~\cite{[HovCBC]}.  The spectrum $\hfp$ is a ring object, hence local with respect to itself.  As mentioned in Section~\ref{sect-harm}, it is known that $\hfp$ is dissonant and $M(p)$ is harmonic, so $[\hfp, M(p)]_n=0$ for all $n$, and $\hfp\in \A$.  On the other hand, the identity on $M(p)$ is nonzero, so $M(p)\notin \A$.  This shows that $\A$ is a localizing subcategory that is not a Bousfield class.

Another example comes from $Z=L_\hfp(BP)$.  Clearly $Z\notin \<Z^*\>$.  But $BP$ is also harmonic, so $[\hfp, BP]_n = 0$ and $[\hfp, Z]_n=0$ for all $n$, and $\hfp\in \<Z^*\>$.  Since $Z\in \<M(p)^*\>$, we know that $\<M(p)^*\> \neq \<Z^*\>$. \end{pf}

Both these counterexamples are cohomological Bousfield classes.  It would be interesting to find a localizing subcategory in $\LL_\hfp$ that is not a cohomological Bousfield class, or show there are none.  Also, it is not clear what, if anything, the previous proposition might tell us about the original conjecture in $\SS$.  For example, as localizing subcategories in $\SS$, we have~\cite[3.3]{[HovCBC]} that $\<M(p)^*\> = \<I\>$.

\subsection{ The $I$-local category}

Recall by $I$ we mean the Brown-Comenetz dual of the sphere spectrum.  It is a rare example of a nonzero spectrum that squares to zero.  Hovey and Palmieri~\cite[Lemma~7.8]{[HP]} conjecture that $\<I\>$ is minimal in $\BL(\SS)$.  

\begin{theorem} \label{i-loc-tc} On $\LL_I$, for all $n$ we have that $\LTCi$, $\LTCii$, and $\LTCiii$ all hold.
\end{theorem}

\begin{pf} By Lemma $7.1(c)$ of~\cite{[HP]}, $T(n)\w I=0$ for all $n$, so $LT(n)=0$.  Since $K(i)$ is a $BP$-module, and $BP\w I=0$ by~\cite[Cor.~B.11]{[HovS]}, we also have $K(n)\w I=0$ for all $n$.  Therefore $\<LT(n)\> = \0=\<LK(n)\>$ for all $n$, and the rest follows from Theorem~\ref{LTCs}. \end{pf}

\begin{proposition} \label{bl-i} The Bousfield lattice of $\LL_I$ is the two-element lattice $\{\0, \<L_IS^0\>\}$.
\end{proposition}

\begin{pf} By~\cite[7.1(c)]{[HP]}, $\<I\><\<\hfp\>$.  Then Proposition~\ref{lattice-map} implies that there is a surjective lattice map from $\BL(\LL_\hfp) = \{\0, \<\hfp\>\}$ onto $\BL(\LL_I)$.  Note that by Lemma~\ref{compose}, $\<L_I\hfp\> = \<L_IL_\hfp S^0\> = \<L_IS^0\>$.  

It remains to show that $\<L_IS^0\>\neq \0$.  But any $X$ in $\SS$ with $X\w I\neq 0$ in $\SS$ will have $L_IX\neq 0$ and $L_IX\notin \<L_IS^0\>$ in $\BL(\LL_I)$; this is because $L_IX\w_{\LL_I}L_IS^0 = L_I(L_IX \w_{\SS} L_IS^0) = L_I(X\w_{\SS} S^0) = L_I(X)$.  For example, $F(n)\w I\neq 0$ for all $n$~\cite[7.1(e)]{[HP]}.   \end{pf}

\begin{corollary} In $\LL_I$, the $\GSC$ fails but the $\SDGSC$ holds. \label{gsc-i}
\end{corollary}

\begin{pf} Corollary B.13 of~\cite{[HovS]} also shows that $\LL_I$ has no nonzero compacts, so the proof is the same as for $\LL_\hfp$. \end{pf}

In~\cite[Conj.~3.10]{[Hov95a]}, Hovey states the {Dichotomy Conjecture}: In $\SS$ every spectrum has either a finite local or a finite acyclic.  In~\cite{[HP]} the authors discuss several equivalent formulations, and some implications.  We briefly point out a relationship between this conjecture and Proposition~\ref{bl-i}.

\begin{proposition} If the Dichotomy Conjecture holds, then the cardinality of $\BL(\LL_I)$ is at most two.
\end{proposition}

\begin{pf} Lemma 7.8 of~\cite{[HP]} shows that if the Dichotomy Conjecture holds, then $\<I\>$ is minimal among nonzero classes in $\BL(\SS)$.  This is the case if and only if $a\<I\>$ is maximal among non-top classes in $\BL(\SS)$, where $a(-)$ denotes the complementation operation first studied by Bousfield~\cite{[Bou79b]}.  Let $a\<I\>\ua$ denote the sublattice $\{\<X\>\;|\; \<X\>\geq a\<I\>\}\subseteq \BL(\SS)$.  In~\cite[Prop.~3.2]{[preprint]} we show that there is a surjective lattice map from $a\<I\>\ua$ onto $\BL(\LL_I)$.  Thus, if the Dichotomy Conjecture holds, $a\<I\>\ua$ has cardinality two and $\BL(\LL_I)$ has cardinality at most two.
\end{pf}

As for classifying localizing subcategories of $\LL_I$, or at least perhaps finding a proper nonzero localizing subcategory, we must get around the fact that so many spectra are $I$-acyclic.  We know that $LF(n)\neq 0$ for all $n$, however $\loc(LF(n))$ is the acyclics of $\l_{n-1}^f:\LL_I\ra \LL_I$ and Theorem~\ref{i-loc-tc} shows that $\lnf=0$ for all $n$.  Thus $\loc(LF(n))=\loc(LS^0)$ in $\LL_I$ for each $n$.

\subsection{The $BP$-local category}

\begin{theorem} \label{LTCs-bp} On $\LL_\bp$, for all $n$ we have that $\LTCi$, $\LTCii$, and $\LTCiii$ all hold. 
\end{theorem}

\begin{pf} We will show that $\LTCii$ holds for all $n$, and the rest follows from Theorem~\ref{LTCs}.  Since each $K(i)$ is a $BP$-module spectrum, $\<K(i)\>\leq \<BP\>$, and since $K(i)$ is local with respect to itself this implies that $K(i)$ is $BP$-local.  Furthermore, this implies $\<E(n)\>\leq \<BP\>$, so from Lemma~\ref{compose} $L_n=L_nL=LL_n$ as functors on $\SS$.  

We claim that $L_n:\LL_\bp\ra \LL_\bp$, taking $LY\mapsto L_nLY=L_nY$, is a smashing localization functor on $\LL_\bp$.  We have $L_n(LY) = L(L_nY) = L(L_nS^0\w_\SS Y) = L(LL_nS^0\w_\SS LY) = L(L_nLS^0\w_\SS LY) = (L_nLS^0)\w_{\LL_\bp} (LY)$.  This shows that on $\LL_\bp$, the localization functor $L_n$ is also given by smashing with the localization of the unit, $L_nLS^0$, and thus is smashing.

Since both $L_n$ and $\ln$ are localization functors on $\LL_\bp$ that annihilate $\<K(0)\vee\cdots\vee K(n)\> = \<LK(0)\vee\cdots\vee LK(n)\>$, they are isomorphic.

On $\SS$, the natural map $\Lnf X\ra \Ln X$ is a $\bp$-equivalence~\cite[Thm.~2.7(iii)]{[RavLaTC]}.  This means that $L\Lnf X = LL_n X$ for all objects $X$ in $\SS$, in particular for all $\bp$-local objects.  Therefore $L\Lnf = L_n = \ln$ is a smashing localization functor on $\LL_\bp$.  The acyclics are $\<L\Lnf(LS^0)\> = \<L\Lnf S^0\> = \<LT(0)\vee\cdots \vee LT(n)\>$.  These are the same acyclics as for $\lnf$, and so we conclude that $\lnf$ and $\ln$ are isomorphic, and the natural map $\lnf X\ra \ln X$ is an isomorphism.
\end{pf}

\begin{proposition} The $\GSC$ fails in $\LL_\bp$. \label{gsc-bp}
\end{proposition}

\begin{pf}  The proof of the last theorem showed that $L_n:\LL_\bp\ra \LL_\bp$ is a (different) smashing localization for each $n$.  However, by~\cite[Cor.~B.13]{[HovS]} the category $\LL_\bp$ has no nonzero compact objects.
\end{pf}

Note that the $\SDGSC$ could still hold, since all the smashing localizations we have identified on $\LL_\bp$ are of the form $L_n = \ln = \lnf$, so are generated by strongly dualizable objects.  The question of finding any other smashing localizations on $\LL_\bp$ is probably at least as hard as doing so on $\SS$, in light of Proposition~\ref{implications}.

All of $\<E(n)\>$, $\<K(n)\>$, $\<\hfp\>$, and $\<I\>$ are ``small" in $\BL(\SS)$, so by Lemma~\ref{compose} it is not surprising that the Bousfield lattices of their localized categories are not very large; this is not true of $\<BP\>$ in $\BL(\SS)$.  We have the following bounds on the Bousfield lattice of the local category.

\begin{proposition} The Bousfield lattice of $\LL_\bp$ has $2^{\aleph_0}\leq |\BL(\LL_\bp)|\leq 2^{2^{\aleph_0}}$. \label{bl-bp}
\end{proposition}

\begin{pf} The second inequality is Corollary~\ref{bl-size}. Since $\<K(i)\>\leq \<\bp\>$ for all $i$, we have $\<Q\> = \<\bigvee_{i\geq 0}K(i)\>\leq \<\bp\>$, and so by Propositions~\ref{lattice-map} and~\ref{jb} we have $|\BL(\LL_\bp)|\geq |\BL(\HH)| =2^{\aleph_0}$.
\end{pf}

\subsection{The $F(n)$-local category}  We conclude with a short discussion of the $F(n)$-local category.

Any smashing localization $L:\SS\ra \SS$ gives a splitting of the Bousfield lattice $$\BL(\SS)\stackrel{\sim}{\longrightarrow} \BL(\LL_{LS^0})\times \BL(\LL_{CS^0}),$$ where $\<X\>\mapsto (\<X\w LS^0\>, \<X\w CS^0\>)$.  See~\cite[Prop.~6.12]{[IK]} or~\cite[Thm.~5.14]{[preprint]} for more details.  If we take $L=\Lnf:\SS\ra\SS$, then we have $\<LS^0\> = \<T(0)\vee\cdots\vee T(n)\>$ and $\<CS^0\> = \<F(n+1)\>$.  Of course, the relationship between $\LL_{T(0)\vee\cdots\vee T(n)}$ and $\LL_{E(n)}$ of Subsection~\ref{sect-en} is immediately related to the original $\TCi$ in $\SS$.  However, this suggests that $\LL_{F(n)}$ is worth investigating further.

By Lemma~\ref{calcs}, in $\BL(\SS)$ there is a chain
$$\<S^0\> = \<F(0)\>\geq \<F(1)\> \geq \<F(2)\>\geq\cdots,$$ and by Lemma~\ref{compose} this gives a chain of lattice surjections
$$\BL(\SS) = \BL(\LL_{F(0)}) \twoheadrightarrow \BL(\LL_{F(1)}) \twoheadrightarrow \BL(\LL_{F(2)})\twoheadrightarrow\cdots.$$

From the above observations, we expect $\BL(\LL_{F(n)})$ to be about as complicated as $\BL(\SS)$.  For example, $F(n)\w I\neq 0$ for all $n$, and so $L_{F(n)}I$ is a square-zero object in $\LL_{F(n)}$.  This means that, unlike in most of the localized categories discussed throughout this paper, we know that $\BA(\LL_{F(n)})\neq \BL(\LL_{F(n)})$.

\end{document}